\documentclass[reqno]{amsart}
\usepackage{amssymb, amsmath, amsthm}
\usepackage{bbold}
\usepackage[colorlinks]{hyperref}
\providecommand{\mathbbm}[1]{\mathbb{#1}}
\allowdisplaybreaks

\theoremstyle{plain}
\newtheorem{theorem}{Theorem}[section]
\theoremstyle{remark}
\newtheorem{remark}[theorem]{Remark}
\theoremstyle{plain}
\newtheorem{lemma}[theorem]{Lemma}
\newtheorem{proposition}[theorem]{Proposition}
\newtheorem{corollary}[theorem]{Corollary}
\newtheorem{definition}[theorem]{Definition}
\newtheorem{hypothesis}[theorem]{Hypothesis}
\numberwithin{equation}{section}
\begin{document}
\raggedbottom
\author{E. M. Ait Benhassi}
\address{CRMEF Marrakech\\
	Lab. of Mathematics, Modeling and Automatic Systems\\
	Marrakech, Morocco}
\email{m.benhassi@uca.ma}

\author{Mohamed Fadili}
\address{Ecole Normale Sup\'erieure, Universit\'e Cadi Ayyad\\
	Lab. of Mathematics, Modeling and Automatic Systems\\
	Marrakech, Morocco}
\email{m.fadili@uca.ma}

\author{Lahcen Maniar}
\address{Faculty of Sciences Semlalia, Cadi Ayyad University\\
	Lab. of Mathematics, Modeling and Automatic Systems\\
	Marrakech, Morocco}
\email{maniar@uca.ma}
\title[Null controllability with variable delay]{Null controllability of degenerate parabolic equations with a time-varying delay}
\thanks{An earlier version of this work treated the constant-delay problem.}

\begin{abstract}
This paper is devoted to the null controllability of a one-dimensional linear
parabolic equation with a time-varying state delay and a diffusion coefficient
that degenerates at one endpoint. Both the weakly and strongly degenerate
cases are considered. When the delay map $g(t)=t-\tau(t)$ is strictly
increasing, the adjoint equation contains an advanced term weighted by the
Jacobian associated with $g^{-1}$. For a distributed control acting on an
interior subinterval, we combine a degenerate Carleman estimate on the final
delay-free interval with a monotonicity functional whose upper integration
limit moves in time. The resulting observability inequality gives null
controllability under a terminal decay condition on the delayed coefficient
outside the control region. We then obtain Dirichlet boundary controllability
at the nondegenerate endpoint by extending the equation in space. The
constant-delay equation is recovered as a special case.
\end{abstract}
\maketitle

\section{Introduction}
We consider the following linear degenerate parabolic equation with a
time-varying delay and an internal localized control:
\begin{equation}\label{sys:internal}
	\begin{cases}
		y_t=\left(a(x)y_x\right)_x+b(t,x)y+c(t,x)y(t-\tau(t),x)+\mathbbm{1}_{\omega}u
		& \text{ in } Q,\\
		y(t,0)=y(t,1)=0 & \text{in the WD case},\\
		(ay_x)(t,0)=y(t,1)=0 & \text{in the SD case},\\
		y(0,\cdot)=y_0 & \text{ in } (0,1),\\
		y=\Theta & \text{ on } (-\tau^*,0)\times(0,1),
	\end{cases}
\end{equation}
The second problem considered in this paper is the boundary controlled system
\begin{equation}\label{sys:boundary}
\begin{cases}
y_t=\left(a(x)y_x\right)_x+b(t,x)y+c(t,x)y(t-\tau(t),x) & \text{in } Q,\\
\begin{cases}
y(t,0)=0, & \text{in the WD case},\\
(ay_x)(t,0)=0, & \text{in the SD case},
\end{cases}\\
y(t,1)=v(t) & \text{for }t\in(0,T),\\
y(0,\cdot)=y_0,\quad y=\Theta \text{ on }(-\tau^*,0)\times(0,1),
\end{cases}
\end{equation}
Here $v\in L^2(0,T)$, $T>0$, $Q=(0,T)\times(0,1)$,
$\omega\Subset(0,1)$ is a nonempty open interval,
$b,c\in L^\infty(Q)$, $y_0\in L^2(0,1)$,
$\Theta\in L^2((-\tau^*,0)\times(0,1))$, and
$u\in L^2((0,T)\times\omega)$. The history $\Theta$ describes the state
before $t=0$.

\begin{hypothesis}[Time-varying delay]\label{hyp:delay}
The function $\tau$ belongs to $W^{1,\infty}(0,T)$ and there exist constants
$\tau_*,\tau^*,\delta_\tau>0$ such that
\[
0<\tau_*\le\tau(t)\le\tau^*,\qquad
1-\tau'(t)\ge\delta_\tau\quad\text{for a.e. }t\in(0,T).
\]
Set $g(t)=t-\tau(t)$. Then $g$ is strictly increasing and its inverse
$\eta=g^{-1}$ is Lipschitz continuous on $[g(0),g(T)]$.
\end{hypothesis}

The diffusion coefficient $a$ is continuous and degenerates at $x=0$, i.e., $a(0)=0$ and $a(x)>0$ for $x\in(0,1]$. The weakly and strongly degenerate regimes, together with their respective boundary conditions, are stated precisely in Section~\ref{sec:degeneracy}.

The approximate controllability of infinite-dimensional retarded linear systems has been investigated in \cite{Nakagiri,Nakagiri&Yamamoto,Curt&Zwart}. More recently, Ammar-Khodja et al.~\cite{A5} established the first null controllability result for retarded nondegenerate parabolic equations with an internally localized control. In the present paper, we follow the strategy of \cite{A5} and combine it with degenerate Carleman estimates \cite{BOU,FadiliManiar,benfama} to treat retarded degenerate parabolic equations. We also study the boundary control problem when the control acts at the nondegenerate endpoint $x=1$.

Degenerate parabolic equations with Volterra memory have subsequently been
studied in \cite{AllalFragnelli,ZhouZhang} and, for degenerate/singular
systems, in \cite{AllalFragnelliSalhi}. Those models contain a time integral
of the past trajectory. The variable discrete delay considered here has a
different adjoint structure. The substitution $r=g(t)$ produces an advanced
term composed with $\eta$ and weighted by $1/g'(\eta(r))$; this term vanishes
on $(\max\{0,g(T)\},T)$. Our argument first applies a degenerate Carleman
estimate on that delay-free interval and then propagates the resulting
information back by a monotonicity functional with a moving upper limit.

Time-varying delays have also been considered in stabilization problems for
parabolic PDEs, for instance through mobile collocated actuators and sensors
in \cite{WuZhang}, and in delayed boundary or internal feedbacks for coupled
systems \cite{Majumdar}. Those works concern asymptotic stabilization. The
present finite-time null controllability problem instead requires an
observability inequality for the nonlocal advanced adjoint equation.

We first derive the exact adjoint and the dual history observation associated
with the nonaffine time map $g$. A moving-limit energy then gives an
observability estimate that is uniform for delays satisfying
Hypothesis~\ref{hyp:delay} with the same structural constants. Finally, we
transfer the distributed result to a Dirichlet control at the nondegenerate
endpoint without requiring additional regularity of the history. To the best
of our knowledge, finite-time controllability has not been studied previously
for a degenerate parabolic equation with a time-varying state delay.

For the boundary problem, the controlled endpoint is nondegenerate because
$a(1)>0$. We extend the equation beyond $x=1$ and place a distributed
control in the added nondegenerate region. The trace of the extended state at
$x=1$ then supplies the Dirichlet control.

The constant-delay problem treated in an earlier version of this work is
recovered when $\tau(t)\equiv h$. The novelty here lies in the inverse delay map, its
Jacobian, and the moving-limit energy, rather than in a new proof of that
special case. As in the nondegenerate case \cite{A5}, reaching zero at an
intermediate time does not generally keep the state at rest. This explains
the decay assumption on $c$ near $T$.

The main results are the following.
\begin{theorem}[Distributed Null Controllability]\label{thm:internal}
Assume Hypothesis~\ref{hyp:WD} or Hypothesis~\ref{hyp:SD} and
Hypothesis~\ref{hyp:delay}, and suppose that
$c$ satisfies
\begin{equation}\label{cond:decay}
\lim_{t\to T^-}(T-t)^4 \ln \|c(t,\cdot)\|_{L^\infty((0,1)\setminus\overline{\omega})} = -\infty.
\end{equation}
Here and below, $\ln0=-\infty$.
Then for every $(y_0,\Theta)\in M_\tau:=L^2(0,1)\times
L^2((-\tau^*,0)\times(0,1))$, there exists
$u\in L^2((0,T)\times\omega)$ such that the solution of
\eqref{sys:internal} satisfies $y(T)=0$ in $(0,1)$. Moreover,
\[
\|u\|_{L^2((0,T)\times\omega)}\le C_T\|(y_0,\Theta)\|_{M_\tau}.
\]
\end{theorem}
\begin{theorem}[Boundary Null Controllability]\label{thm:boundary}
Assume Hypothesis~\ref{hyp:WD} or Hypothesis~\ref{hyp:SD} and
Hypothesis~\ref{hyp:delay}. Suppose that $c$
satisfies the global decay condition
\begin{equation}\label{cond:decay-bd}
\lim_{t\to T^-}(T-t)^4 \ln \|c(t,\cdot)\|_{L^\infty(0,1)} = -\infty.
\end{equation}
Then for every $(y_0,\Theta)\in M_\tau$, there exists $v\in L^2(0,T)$ such that the solution of \eqref{sys:boundary} satisfies $y(T)=0$ in $(0,1)$. Moreover,
\[
\|v\|_{L^2(0,T)}\le C_T\|(y_0,\Theta)\|_{M_\tau}.
\]
\end{theorem}
In both theorems, the constant $C_T$ is independent of $(y_0,\Theta)$ and of
the particular delay function $\tau$; it may depend on the control region and
on $a,b,c,T,\tau_*,\tau^*$ and $\delta_\tau$.
\begin{corollary}[Constant delay]\label{cor:constant-delay}
If $\tau(t)\equiv h>0$, then $g(t)=t-h$, $\eta(s)=s+h$, and
$\delta_\tau=1$. Theorems~\ref{thm:internal} and~\ref{thm:boundary} reduce
to the corresponding constant-delay controllability results.
\end{corollary}
\begin{proof}
In this case $T_\star=\max\{0,T-h\}$ and
$q_W(t)=(cW\mathbbm1_{[0,T]})(t+h)$. Moreover,
$\beta(t)=\min\{t+h,T\}$, so the moving-limit energy below becomes the
constant-delay energy. The state equation, the history term, and both decay
conditions therefore coincide with their constant-delay counterparts.
\end{proof}
\begin{remark}[Nonconstant admissible delays]\label{rem:examples}
For example, let
\[
\tau(t)=h+\varepsilon\sin(2\pi t/T),
\]
where $h>\varepsilon>0$ and $2\pi\varepsilon/T<1$. Then
Hypothesis~\ref{hyp:delay} holds with
$\tau_*=h-\varepsilon$, $\tau^*=h+\varepsilon$, and
$\delta_\tau=1-2\pi\varepsilon/T$. The terminal condition is also nonempty.
For instance, outside the control region one may take
\[
c(t,x)=d(t,x)\exp\!\left(-\frac{1}{(T-t)^p}\right),
\qquad p>4,
\]
with $d\in L^\infty(Q)$. Inside $\omega$, the coefficient may remain merely
bounded. For boundary controllability, the displayed decay is imposed on the
whole spatial interval.
\end{remark}
Section 2 states the degeneracy hypotheses and fixes the functional setting.
In Section 3, we prove well-posedness, derive the adjoint formulation and the
Carleman estimate, absorb the delay term, and establish distributed null
controllability. Section 4 treats boundary null controllability by extending
the equation beyond the nondegenerate endpoint. We close with a brief
conclusion.

Throughout the paper, generic positive constants may change from line to line. We write $\|\cdot\|$ for the $L^2(0,1)$ norm and $\|\cdot\|_X$ for the norm in a Banach space $X$.

\section{Degeneracy hypotheses and notation}\label{sec:degeneracy}
We use one of the following two assumptions on the diffusion coefficient.

\begin{hypothesis}[Weak degeneracy (WD)]\label{hyp:WD}
The coefficient satisfies $a\in C([0,1])\cap C^1((0,1])$, $a(0)=0$, and $a>0$ on $(0,1]$. Moreover, there exists $K_a\in[0,1)$ such that
\[
xa'(x)\le K_a a(x),\qquad x\in(0,1].
\]
The associated boundary conditions are $y(t,0)=y(t,1)=0$.
\end{hypothesis}

\begin{hypothesis}[Strong degeneracy (SD)]\label{hyp:SD}
The coefficient satisfies $a\in C^1([0,1])$, $a(0)=0$, and $a>0$ on $(0,1]$. Moreover, there exists $K_a\in[1,2)$ such that
\[
xa'(x)\le K_a a(x),\qquad x\in[0,1].
\]
If $K_a>1$, there exists $\vartheta\in(1,K_a]$ such that $x\mapsto a(x)/x^\vartheta$ is nondecreasing near $0$; if $K_a=1$, the same property holds for some $\vartheta\in(0,1)$. The associated boundary conditions are $(ay_x)(t,0)=y(t,1)=0$.
\end{hypothesis}

The energy space depends on the degeneracy regime. In the weakly degenerate
case, set
\[
H_{a,\mathrm{WD}}^1(0,1)=
\left\{z\in AC([0,1]):\sqrt a\,z_x\in L^2(0,1),\quad
z(0)=z(1)=0\right\}.
\]
In the strongly degenerate case, set
\[
H_{a,\mathrm{SD}}^1(0,1)=
\left\{z\in L^2(0,1):z\in AC_{\mathrm{loc}}((0,1]),\quad
\sqrt a\,z_x\in L^2(0,1),\quad z(1)=0\right\}.
\]
Both spaces are endowed with the norm
\[
\|z\|_{H_a^1}^2=\|z\|_{L^2(0,1)}^2+
\|\sqrt a\,z_x\|_{L^2(0,1)}^2.
\]
We write $H_a^1$ for the relevant one of these spaces and define
\[
D(A)=\left\{z\in H_a^1:az_x\in H^1(0,1)\right\},
\qquad Az=(az_x)_x,
\]
with $(az_x)(0)=0$ imposed in the strongly degenerate case. We also write
$H_a^2=D(A)$ and denote by $H_a^{-1}$ the dual of $H_a^1$ with pivot
$L^2(0,1)$.

\begin{definition}[Temporal weight]\label{def:temporal-weight}
Set $T_\star=\max\{0,T-\tau_*\}$. Since
$g(T)=T-\tau(T)\le T-\tau_*$, the advanced term in the adjoint equation
below vanishes on $(T_\star,T)$. Define
\[
\theta(t)=\frac{1}{(t-T_\star)^4(T-t)^4},\qquad t\in(T_\star,T).
\]
The reciprocal weight $\theta^{-1}$ vanishes at $T_\star$ and $T$, while
$\theta$ reaches its minimum at $(T+T_\star)/2$. Moreover,
$T-T_\star=\min\{T,\tau_*\}>0$.
\end{definition}

\section{Distributed control}
\subsection{Well-posedness}
The retarded evolution theory of Artola \cite{Artola}, together with the
arguments in \cite{A5}, gives the following well-posedness result.
\begin{proposition}\label{prop:well-int}
If $(y_0,\Theta,u)\in M_\tau\times L^2(Q)$, then
\eqref{sys:internal} admits a unique mild solution
$y\in L^2((-\tau^*,T)\times(0,1))\cap L^2(0,T;H_a^1)
\cap \mathcal{C}([0,T];L^2(0,1))$, with
$y_t\in L^2(0,T;H_a^{-1})$, satisfying the energy estimate
\begin{align*}
&\sup_{t\in[0,T]}\|y(t)\|_{L^2}^2
 + \int_0^T\left( \|\sqrt{a}\,y_x\|_{L^2}^2 + \|y_t\|_{H_a^{-1}}^2 \right)dt \\
&\qquad\le C_T\left( \|y_0\|_{L^2}^2 + \|\Theta\|_{L^2}^2 + \|u\|_{L^2(Q)}^2 \right).
\end{align*}
Furthermore, if $(y_0,\Theta,u)\in H_a^1(0,1)\times
L^2((-\tau^*,0)\times(0,1))\times L^2(Q)$, then
\[
y\in L^2(0,T;H_a^2)\cap \mathcal{C}([0,T];H_a^1(0,1)),\qquad
y_t\in L^2(0,T;L^2(0,1)),
\]
and
\begin{align*}
&\sup_{t\in[0,T]}\|y(t)\|_{H_a^1}^2
+\int_0^T\left(\|y_t\|_{L^2}^2+\|(ay_x)_x\|_{L^2}^2\right)dt\\
&\qquad\le C_T\left(\|y_0\|_{H_a^1}^2+
\|\Theta\|_{L^2((-\tau^*,0)\times(0,1))}^2+
\|u\|_{L^2(Q)}^2\right),
\end{align*}
where $C_T>0$ is independent of $(y_0,\Theta,u)$.
\end{proposition}
\begin{proof}
Let $A$ be the operator defined on $L^2(0,1)$ by $Ay=(ay_x)_x$ with domain
\[
D(A) = \{y\in H_a^1(0,1)\mid (ay_x)_x\in L^2(0,1)\},
\]
where the boundary conditions are those of Hypothesis~\ref{hyp:WD} or Hypothesis~\ref{hyp:SD}. Under either degeneracy assumption, integration by parts shows that $A$ is self-adjoint and dissipative and that its closed form is
$\mathfrak a(z,z)=\int_0^1a|z_x|^2dx$ on $H_a^1(0,1)$ with the corresponding boundary conditions. Hence $A$ generates a contraction analytic semigroup on $L^2(0,1)$.

Let $I_k=[k\tau_*,\min\{(k+1)\tau_*,T\}]$. For $t\in I_k$,
\[
g(t)=t-\tau(t)\le t-\tau_*\le k\tau_*.
\]
Thus $y(g(t))$ is either prescribed by the history or belongs to a portion of
the trajectory already constructed. On each $I_k$ the equation is therefore
a nonautonomous parabolic equation without delay and with a known forcing.
Moreover, the change of variables $r=g(t)$ and $g'\ge\delta_\tau$ yield
\[
\int_J\|z(g(t))\|^2dt
\le\delta_\tau^{-1}\int_{g(J)}\|z(r)\|^2dr
\]
for every interval $J$ and every admissible history or previously constructed
trajectory $z$. Hence the forcing belongs to $L^2$. The variation-of-constants
theorem applies successively on at most $\lceil T/\tau_*\rceil$ intervals,
giving existence and uniqueness.

For the estimates, multiply the equation by $y$ and use the boundary conditions. Young's inequality gives, for almost every $t$,
\begin{align*}
\frac{d}{dt}\|y(t)\|^2+\|\sqrt{a}\,y_x(t)\|^2
&\le C\|y(t)\|^2+C\|y(g(t))\|^2
+C\|u(t)\|_{L^2(\omega)}^2.
\end{align*}
The preceding change of variables also gives
\[
\int_0^T\|y(g(t))\|^2dt
\le\delta_\tau^{-1}\left(
\|\Theta\|_{L^2((-\tau^*,0)\times(0,1))}^2+
\int_0^T\|y(t)\|^2dt\right).
\]
Integration of the differential inequality and Gronwall's lemma therefore yield
\[
\sup_{0\le t\le T}\|y(t)\|^2
+\int_0^T\|\sqrt{a}\,y_x(t)\|^2dt
\le C_T\bigl(\|y_0\|^2+\|\Theta\|^2+\|u\|^2\bigr).
\]
The equation, viewed in $H_a^{-1}$, then gives the asserted bound for $y_t$.

If $y_0\in H_a^1=D((-A)^{1/2})$, maximal $L^2$ regularity for the analytic semigroup gives on $I_0$
\begin{align*}
&\sup_{t\in I_0}\|y(t)\|_{H_a^1}^2
+\int_{I_0}\bigl(\|y_t\|^2+\|Ay\|^2\bigr)dt\\
&\qquad\le C_T\bigl(\|y_0\|_{H_a^1}^2+
\|c(\cdot)y(g(\cdot))+\mathbbm1_\omega u\|_{L^2(I_0;L^2)}^2\bigr).
\end{align*}
Repeating this estimate on the intervals $I_k$ and using the preceding energy bound for the delayed forcing gives the second estimate. Since $Ay=(ay_x)_x$ and $D((-A)^{1/2})=H_a^1$, the stated regularity follows.
\end{proof}

\subsection{Adjoint formulation and observability criterion}
For a function $W$ on $Q$, define
\[
q_W(t,x)=
\begin{cases}
\displaystyle
\frac{c(\eta(t),x)W(\eta(t),x)}{g'(\eta(t))},
&0\le t\le g(T),\\[2mm]
0,&\max\{0,g(T)\}<t\le T,
\end{cases}
\]
where the first branch is empty when $g(T)<0$. The adjoint equation
associated with the distributed system is
\begin{equation}\label{eq:adjoint}
\begin{cases}
-W_t = (a(x)W_x)_x + b(t,x)W + q_W(t,x) & \text{in } Q,\\
W(t,0)=W(t,1)=0 & \text{in the WD case},\\
(aW_x)(t,0)=W(t,1)=0 & \text{in the SD case},\\
W(T,\cdot)=W_0 & \text{in } (0,1).
\end{cases}
\end{equation}
This problem is well posed backward in time. On the final interval
$(T_\star,T)$ one has $q_W=0$. On every preceding interval of length
$\tau_*$, the value $\eta(t)$ satisfies
\(\eta(t)=t+\tau(\eta(t))\ge t+\tau_*\), so it belongs to a portion of the
trajectory already constructed while proceeding from $T$ to $0$. Moreover,
if $g(T)\ge0$, the change of variables $s=g(r)$ gives
\[
\int_0^{\max\{0,g(T)\}}\|q_W(s)\|^2ds
=\int_{\eta(0)}^T\frac{\|c(r)W(r)\|^2}{g'(r)}dr
\le\frac{\|c\|_\infty^2}{\delta_\tau}
\int_0^T\|W(r)\|^2dr.
\]
If $g(T)<0$, then $q_W=0$ on all of $(0,T)$.
The same parabolic estimates as in Proposition~\ref{prop:well-int} therefore
give a unique adjoint solution in the corresponding energy class.
Let $Q_\omega=(0,T)\times\omega$. We define the bounded operators
\[
S_T:M_\tau\longrightarrow L^2(0,1),\qquad
S_T(y_0,\Theta)=y^H(T),
\]
where $y^H$ is the uncontrolled solution, and
\[
L_T:L^2(Q_\omega)\longrightarrow L^2(0,1),\qquad
L_Tu=y^u(T),
\]
where $y^u$ has zero initial state and zero history. By linearity,
\[
y(T)=S_T(y_0,\Theta)+L_Tu.
\]

\begin{proposition}[Observability criterion]\label{prop:obs-criterion}
Define, for $s\in(-\tau^*,0)$,
\[
H_W(s,x)=
\begin{cases}
\displaystyle
\frac{c(\eta(s),x)W(\eta(s),x)}{g'(\eta(s))},
&g(0)\le s\le\min\{0,g(T)\},\\[2mm]
0,&\text{otherwise}.
\end{cases}
\]
The history observation is well defined and depends continuously on $W$.
Indeed, changing variables $s=g(r)$ gives
\[
\|H_W\|^2
=\int_0^{\beta(0)}\frac{\|c(r)W(r)\|^2}{g'(r)}dr
\le\frac{\|c\|_\infty^2}{\delta_\tau}
\|W\|_{L^2(0,T;L^2(0,1))}^2.
\]
System \eqref{sys:internal} is null controllable at time $T$ if and only if
there exists $C_T>0$ such that every solution of \eqref{eq:adjoint} satisfies
\begin{equation}\label{eq:obs-int}
\|W(0)\|_{L^2(0,1)}^2+
\|H_W\|_{L^2((-\tau^*,0)\times(0,1))}^2
\le C_T\int_0^T\int_\omega W^2dxdt,
\end{equation}
for every $W_0\in L^2(0,1)$.
\end{proposition}
\begin{proof}
Let $y$ solve \eqref{sys:internal} and let $W$ solve \eqref{eq:adjoint}. We first take smooth data. Multiplication of the state equation by $W$, followed by integration over $Q$, gives
\begin{align}
\int_Q y_tW\,dx\,dt
&=\int_Q (ay_x)_xW\,dx\,dt+\int_Q byW\,dx\,dt \notag\\
&\quad+\int_Q c(t,x)y(g(t),x)W(t,x)\,dx\,dt
+\int_{Q_\omega}uW\,dx\,dt.\label{eq:duality-start}
\end{align}
Integration by parts in time and space is legitimate for smooth data; the
boundary terms vanish under the homogeneous conditions of
Hypothesis~\ref{hyp:WD} or Hypothesis~\ref{hyp:SD}. Hence
\begin{align*}
\int_Qy_tW\,dx\,dt
&=\langle y(T),W_0\rangle-\langle y_0,W(0)\rangle-\int_QyW_t\,dx\,dt,\\
\int_Q(ay_x)_xW\,dx\,dt&=\int_Qy(aW_x)_x\,dx\,dt.
\end{align*}
For the delayed integral, set $s=g(t)$. Since $g$ is strictly increasing,
the one-dimensional change-of-variables formula gives
\begin{align*}
&\int_0^T c(t,x)W(t,x)y(g(t),x)\,dt\\
&=\int_{g(0)}^{g(T)}y(s,x)
\frac{c(\eta(s),x)W(\eta(s),x)}{g'(\eta(s))}\,ds\\
&=\int_{-\tau^*}^0\Theta(s,x)H_W(s,x)\,ds
+\int_0^Ty(s,x)q_W(s,x)\,ds.
\end{align*}
Substitution in \eqref{eq:duality-start} and use of \eqref{eq:adjoint} cancel the last integral. We obtain the exact duality identity
\begin{align}
\langle y(T),W_0\rangle
&=\langle y_0,W(0)\rangle
+\int_{-\tau^*}^0\!\!\int_0^1\Theta(s,x)H_W(s,x)\,dx\,ds \notag\\
&\quad+\int_{Q_\omega}uW\,dx\,dt.\label{eq:duality-int}
\end{align}
Density and the estimates of Proposition \ref{prop:well-int} extend \eqref{eq:duality-int} to all admissible data. Taking successively $u=0$ and $(y_0,\Theta)=(0,0)$ identifies the adjoints:
\[
S_T^*W_0=\bigl(W(0),H_W\bigr),
\qquad L_T^*W_0=W|_{Q_\omega}.
\]
Null controllability is equivalent to $\operatorname{Ran}S_T\subset\operatorname{Ran}L_T$. The Hilbert-space range inclusion theorem \cite{Zab} states that this is equivalent to
\[
\|S_T^*W_0\|_{M_\tau}^2\le C_T\|L_T^*W_0\|_{L^2(Q_\omega)}^2
\quad\text{for all }W_0\in L^2(0,1),
\]
which is precisely \eqref{eq:obs-int}.
\end{proof}

\subsection{Carleman estimates and observability}
With $T_\star$ and $\theta$ as in Definition~\ref{def:temporal-weight}, define
\[
 \psi(x)=\lambda\left(\int_0^x\frac{r}{a(r)}\,dr-d\right),
 \qquad \varphi(t,x)=\theta(t)\psi(x),
\]
where $d>4\int_0^1r/a(r)\,dr$. Choose an interval
$\omega_0\Subset\omega$ and a function $\sigma\in C^2([0,1])$ such that
$\sigma>0$ in $(0,1)$, $\sigma(0)=\sigma(1)=0$, and
$\sigma_x\ne0$ on $[0,1]\setminus\omega_0$. For suitable positive
parameters $\rho$ and $\lambda$, set
\[
 \Psi(x)=e^{\rho\sigma(x)}-e^{2\rho\|\sigma\|_\infty},
 \qquad \Phi(t,x)=\theta(t)\Psi(x).
\]
In particular, $\psi<0$ and $\Psi<0$ on $[0,1]$. The parameters can be
chosen so that the standard comparison conditions between the two weights
hold; see \cite[Propositions~3.4 and~3.5]{FadiliManiar} and
\cite[Section~3.1]{benfama}. We use the following scalar consequence.

\begin{theorem}[Carleman estimates]\label{thm:carleman}
Let $\omega\Subset(0,1)$ and $b\in L^\infty((T_\star,T)\times(0,1))$.
There exist $C,s_0>0$ such that every sufficiently regular solution of
\[
 z_t-(az_x)_x+bz=f
\]
satisfying the boundary conditions of Hypothesis~\ref{hyp:WD} or
Hypothesis~\ref{hyp:SD} satisfies
\begin{multline}
\int_{T_\star}^T\!\!\int_0^1 \big(s\theta a z_x^2 + (s\theta)^3 \tfrac{x^2}{a}z^2\big)e^{2s\varphi}dxdt \\
\le C\Big( \int_{T_\star}^T\!\!\int_0^1 f^2 e^{2s\Phi}dxdt + \int_{T_\star}^T\!\!\int_\omega (s\theta)^3 z^2 e^{2s\Phi}dxdt \Big)
\end{multline}
for all $s\ge s_0$. By density, the estimate extends to the natural weak
solution class.
\end{theorem}

\subsection{Observability for the distributed control}
Applying the internal estimate of Theorem \ref{thm:carleman} with $f=0$ yields:
\begin{lemma}\label{lem:obs-local-int}
For $s$ large enough,
\[
\int_{T_\star}^T e^{-2sM\theta} \int_0^1 W^2 dxdt \le C s^3 \int_{T_\star}^T\!\!\int_\omega W^2 dxdt,
\]
where $M=\max_{x\in[0,1]}|\psi(x)|$.
\end{lemma}
\begin{proof}
On $(T_\star,T)\times(0,1)$ one has $t>g(T)$, so $q_W(t)=0$ and
\[
W_t+(aW_x)_x+bW=0.
\]
After the time reversal
$z(t,x)=W(T+T_\star-t,x)$, the function $z$ satisfies the equation in
Theorem~\ref{thm:carleman} with the bounded coefficient
$-b(T+T_\star-t,x)$. The theorem therefore gives
\[
\int_{T_\star}^T\!\!\int_0^1
\left(s\theta aW_x^2+(s\theta)^3\frac{x^2}{a}W^2\right)e^{2s\varphi}dxdt
\le C\int_{T_\star}^T\!\!\int_\omega (s\theta)^3 W^2e^{2s\Phi}dxdt .
\]
To recover the unweighted $L^2$ norm, write pointwise in $x$
\[
W^2=\left(\frac{a^{1/3}}{x^{2/3}}W^2\right)^{3/4}
\left(\frac{x^2}{a}W^2\right)^{1/4}.
\]
Young's inequality implies
\[
\int_0^1W^2dx\le C\int_0^1\frac{a^{1/3}}{x^{2/3}}W^2dx
+C\int_0^1\frac{x^2}{a}W^2dx.
\]
Set $p(x)=x^{4/3}a(x)^{1/3}$. The degeneracy assumptions imply
$p\le Ca$: indeed, $x^2/a(x)$ is bounded in the WD case and is
nondecreasing near zero in the SD case, hence bounded on $(0,1)$. Moreover,
where $a$ is differentiable,
\[
 \frac{xp'(x)}{p(x)}=\frac43+\frac13\frac{xa'(x)}{a(x)}
 \le \frac{4+K_a}{3}<2.
\]
The weighted Hardy--Poincar\'e inequality
\cite[Proposition~2.1]{BOU}, applied with the corresponding endpoint
condition, consequently gives
\[
\int_0^1\frac{a^{1/3}}{x^{2/3}}W^2dx
=\int_0^1\frac{p(x)}{x^2}W^2dx
\le C\int_0^1aW_x^2dx.
\]
After enlarging $s_0$ so that $s\theta\ge1$, we obtain
\[
\int_0^1W^2dx\le C\int_0^1
\left(s\theta aW_x^2+(s\theta)^3\frac{x^2}{a}W^2\right)dx.
\]
Since $\varphi(t,x)\ge-M\theta(t)$,
\begin{align*}
\int_{T_\star}^Te^{-2sM\theta}\int_0^1W^2dxdt
&\le C\int_{T_\star}^T\!\!\int_0^1
\left(s\theta aW_x^2+(s\theta)^3\frac{x^2}{a}W^2\right)e^{2s\varphi}dxdt\\
&\le C\int_{T_\star}^T\!\!\int_\omega(s\theta)^3W^2e^{2s\Phi}dxdt.
\end{align*}
Finally, $\Phi<0$ and $\theta^3e^{2s\Phi}$ tends to zero at both endpoints of $(T_\star,T)$; hence it is bounded there. This gives the stated estimate.
\end{proof}

\subsection{Delay absorption via monotonicity}
We absorb the delay term through the following monotonicity functional.
\begin{lemma}[Energy monotonicity]\label{lem:energy-mono}
Set
\[
\beta(t)=
\begin{cases}
\eta(t),&0\le t\le g(T),\\
T,&g(T)<t\le T,
\end{cases}
\]
where the first branch is empty if $g(T)<0$, and let
$\kappa := 1 + 2\|b\|_\infty +
\delta_\tau^{-1}\|c\|_\infty^2$. Then the energy functional
\[
E(t)=e^{\kappa t}\left(\|W(t)\|^2+
\int_t^{\beta(t)}\frac{\|c(s)W(s)\|^2}{g'(s)}\,ds\right)
\]
is nondecreasing on $[0,T]$.
\end{lemma}
\begin{proof}
We first take $W_0\in H_a^1(0,1)$, set $\psi_d=cW$, and let
\[
E_1(t)=\|W(t)\|^2+
\int_t^{\beta(t)}\frac{\|\psi_d(s)\|^2}{g'(s)}\,ds.
\]
For $t<g(T)$, the inverse function formula gives
$\beta'(t)=1/g'(\eta(t))$ and hence
\[
\frac{d}{dt}\int_t^{\eta(t)}
\frac{\|\psi_d(s)\|^2}{g'(s)}ds
=\|q_W(t)\|^2-\frac{\|\psi_d(t)\|^2}{g'(t)}.
\]
For $t>g(T)$, the same derivative equals
$-\|\psi_d(t)\|^2/g'(t)$ and $q_W(t)=0$. The adjoint equation and integration
by parts therefore give, almost everywhere,
\begin{align*}
E_1'(t)
&=2\int_0^1aW_x^2\,dx-2\int_0^1bW^2\,dx
 -2\langle W,q_W\rangle+\|q_W\|^2
 -\frac{\|c(t)W(t)\|^2}{g'(t)}.
\end{align*}
Completing the square,
\[
-2\langle W,q_W\rangle+\|q_W\|^2
=\|q_W-W\|^2-\|W\|^2
\]
whenever $q_W$ is active; when it vanishes, the corresponding pair is zero
and the same lower bound remains true. Since $g'\ge\delta_\tau$, it follows
that
\[
E_1'(t)\ge-\kappa\|W(t)\|^2.
\]
Since the integral term in $E_1$ is nonnegative,
\[
\kappa E_1(t)+E_1'(t)\ge0.
\]
Thus $E'(t)=e^{\kappa t}(\kappa E_1(t)+E_1'(t))\ge0$. For $W_0\in L^2(0,1)$, approximate by terminal data in $H_a^1(0,1)$. Continuity of the adjoint solution in $C([0,T];L^2(0,1))$ and dominated convergence in the history term pass the monotonicity inequality to the limit.
\end{proof}

\begin{lemma}[Global observability]\label{lem:global-obs}
Under \eqref{cond:decay}, every solution of \eqref{eq:adjoint} satisfies \eqref{eq:obs-int}.
\end{lemma}
\begin{proof}
Put $\ell=T-T_\star>0$, $\nu=\ell/4$, $a_0=T_\star+\nu$, and fix a Carleman parameter $\mathfrak{s}>s_0$.

\medskip
\noindent\textbf{Step A: Local observability estimate.}
Lemma \ref{lem:obs-local-int} gives
\begin{equation}\label{eq:local-used}
\int_{a_0}^{T}e^{-2\mathfrak{s}M\theta(t)}\|W(t)\|^2dt
\le C\mathfrak{s}^3\int_0^T\!\!\int_\omega W^2dxdt.
\end{equation}
Since $e^{-\kappa t}E(t)$ equals the sum of $\|W(t)\|^2$ and the future
integral in the definition of $E$, we have
\begin{align}
\int_{a_0}^{T}e^{-2\mathfrak{s}M\theta(t)-\kappa t}E(t)dt
&\le C\mathfrak{s}^3\int_0^T\!\!\int_\omega W^2dxdt+I,\label{eq:energy-upper}\\
I&:=\int_{a_0}^{T}e^{-2\mathfrak{s}M\theta(t)}
\int_t^{\beta(t)}\frac{\|c(\sigma)W(\sigma)\|^2}{g'(\sigma)}d\sigma\,dt.
\end{align}

\medskip
\noindent\textbf{Step B: Estimating the delay integral $I$.}
Condition \eqref{cond:decay} is equivalent to the following statement: for every $r>0$ there are $\delta>0$ and $C_r>0$ such that
\begin{equation}\label{eq:c-exp-decay}
\|c(t,\cdot)\|_{L^\infty((0,1)\setminus\overline\omega)}^2
\le C_r\exp\left(-\frac{2r}{(T-t)^4}\right),
\qquad T_\star<t<T.
\end{equation}
Indeed, the limit gives this estimate on $(T-\delta,T)$; on $(T_\star,T-\delta)$ it follows after increasing $C_r$, since $c\in L^\infty(Q)$.

Because the integration domain in $I$ is contained in
$\{(t,\sigma):a_0<t<\sigma<T\}$ and $T-a_0=3\nu$, Tonelli's theorem and $e^{-2\mathfrak{s}M\theta}\le1$ imply
\begin{align*}
I&\le\frac{3\nu}{\delta_\tau}
\int_{a_0}^{T}\|c(\sigma)W(\sigma)\|^2d\sigma\\
&\le\frac{3\nu C_r}{\delta_\tau}\int_{a_0}^{T}
e^{-2r/(T-\sigma)^4}\|W(\sigma)\|^2d\sigma
+\frac{3\nu}{\delta_\tau}\|c\|_\infty^2
\int_0^T\!\!\int_\omega W^2dxdt.
\end{align*}

\medskip
\noindent\textbf{Step C: Absorbing the weight.}
Choose
\[
r\ge\frac{\mathfrak{s}M}{\nu^4}.
\]
For $\sigma\in[a_0,T)$ one has $\sigma-T_\star\ge\nu$, and therefore
\[
e^{-2r/(T-\sigma)^4}
\le e^{-2\mathfrak{s}M/((\sigma-T_\star)^4(T-\sigma)^4)}
=e^{-2\mathfrak{s}M\theta(\sigma)}.
\]
Using \eqref{eq:local-used}, we conclude that
\begin{equation}\label{eq:delay-bound}
I\le\frac{3\nu}{\delta_\tau}
\left(C_rC\mathfrak{s}^3+\|c\|_\infty^2\right)
\int_0^T\!\!\int_\omega W^2dxdt.
\end{equation}
Combining \eqref{eq:energy-upper} and \eqref{eq:delay-bound} gives
\begin{equation}\label{eq:weighted-energy-bound}
\int_{a_0}^{T}e^{-2\mathfrak{s}M\theta(t)-\kappa t}E(t)dt
\le C_0\int_0^T\!\!\int_\omega W^2dxdt,
\end{equation}
where $C_0$ depends only on the structural data and $T$.

\medskip
\noindent\textbf{Step D: Truncation and conclusion.}
By Lemma \ref{lem:energy-mono}, $E(0)\le E(t)$. On
$[T_\star+\nu,T-\nu]$ one has
\[
\theta(t)\le \frac{4^8}{3^4\ell^8}.
\]
The length of this interval is $\ell/2$; hence \eqref{eq:weighted-energy-bound} yields
\[
\frac{\ell}{2}\exp\left(-\frac{2^{17}\mathfrak{s}M}{3^4\ell^8}-\kappa T\right)E(0)
\le C_0\int_0^T\!\!\int_\omega W^2dxdt.
\]
Finally, the change of variables $s=g(\sigma)$ gives
\begin{align*}
\|H_W\|_{L^2((-\tau^*,0)\times(0,1))}^2
&=\int_{g(0)}^{\min\{0,g(T)\}}
\left\|\frac{c(\eta(s))W(\eta(s))}{g'(\eta(s))}\right\|^2ds\\
&=\int_0^{\beta(0)}\frac{\|c(\sigma)W(\sigma)\|^2}{g'(\sigma)}d\sigma.
\end{align*}
Consequently,
\[
E(0)=\|W(0)\|^2+\|H_W\|^2,
\]
which is exactly the left-hand side of \eqref{eq:obs-int}.
\end{proof}

\subsection{Proof of distributed null controllability}
\begin{proof}[Proof of Theorem \ref{thm:internal}]
Lemma \ref{lem:global-obs} gives \eqref{eq:obs-int}. By Proposition \ref{prop:obs-criterion}, this inequality is equivalent to
$\operatorname{Ran}S_T\subset\operatorname{Ran}L_T$ and hence to null controllability. More explicitly, the quantitative form of the range inclusion theorem gives a bounded operator
$G_T:M_\tau\to L^2(Q_\omega)$ such that
\[
L_TG_T=-S_T,
\qquad \|G_T(y_0,\Theta)\|_{L^2(Q_\omega)}
\le C_T\|(y_0,\Theta)\|_{M_\tau}.
\]
For prescribed $(y_0,\Theta)$, choose $u=G_T(y_0,\Theta)$ and extend it by zero outside $Q_\omega$. Then
\[
y(T)=S_T(y_0,\Theta)+L_Tu=0,
\]
and the stated estimate for the control follows from the bound on $G_T$.
\end{proof}

\section{Boundary control by extension}\label{sec:boundary}

We prove Theorem~\ref{thm:boundary} by reducing the boundary control problem
to a distributed control problem on a larger interval. This reduction uses
the observability result already obtained. An interior trace estimate also
provides the required \(L^2\)-regularity of the boundary control.

Choose a nonempty open interval \(\omega_e\Subset(1,2)\). We first extend
\(a\) to a coefficient
\(\widetilde a\in C([0,2])\cap C^1((0,2])\) such that
\(\widetilde a=a\) on \([0,1]\), \(\widetilde a>0\) on \((0,2]\),
and \(\widetilde a\) is constant near \(x=2\). The extension is
constructed in Lemma~\ref{lem:extension}. Set
\[
 \widetilde b(t,x)=
 \begin{cases}
  b(t,x),&0<x<1,\\
  0,&1<x<2,
 \end{cases}
\]
and define \(\widetilde c\), \(\widetilde y_0\), and
\(\widetilde\Theta\) by extending \(c\), \(y_0\), and \(\Theta\) by
zero on \((1,2)\). Consider
\begin{equation}\label{sys:extended}
\begin{cases}
 \widetilde y_t=(\widetilde a\widetilde y_x)_x+
 \widetilde b\widetilde y+
 \widetilde c\,\widetilde y(t-\tau(t))+
 \mathbbm{1}_{\omega_e}\widetilde u
 &\text{in }(0,T)\times(0,2),\\
 \widetilde y(t,0)=\widetilde y(t,2)=0
 &\text{in the WD case},\\
 (\widetilde a\widetilde y_x)(t,0)=\widetilde y(t,2)=0
 &\text{in the SD case},\\
 \widetilde y(0,\cdot)=\widetilde y_0,\qquad
 \widetilde y=\widetilde\Theta\ \text{on }(-\tau^*,0)\times(0,2).
\end{cases}
\end{equation}

\begin{lemma}\label{lem:extension}
The coefficient \(a\) admits an extension \(\widetilde a\) with the
properties stated above which satisfies the same type of degeneracy hypothesis
as \(a\). After an inessential rescaling of the spatial interval, the
extended problem therefore falls within the framework of
Theorem~\ref{thm:internal}. Moreover, condition~\eqref{cond:decay-bd}
implies
\[
 \lim_{t\to T^-}(T-t)^4
 \log\|\widetilde c(t,\cdot)\|_
 {L^\infty((0,2)\setminus\overline{\omega_e})}=-\infty.
\]
\end{lemma}

\begin{proof}
Put \(q=a'(1)/a(1)\), so that \(q\le K_a\). Choose
\(\widetilde K_a\) such that
\[
 K_a<\widetilde K_a<1\quad\text{in the WD case},
 \qquad
 K_a<\widetilde K_a<2\quad\text{in the SD case}.
\]
Choose \(\delta>0\) so small that
\((1+\delta)q<\widetilde K_a\) when \(q>0\), and take a smooth
nondecreasing function \(\rho\) on \([0,1]\) such that
\(\rho(0)=0\), \(\rho'(0)=1\), \(0\le\rho'\le1\), and
\(\rho'=0\) on \([\delta,1]\). Define
\[
 \widetilde a(x)=
 \begin{cases}
  a(x),&0\le x\le1,\\
  a(1)\exp\!\bigl(q\rho(x-1)\bigr),&1<x\le2.
 \end{cases}
\]
Then \(\widetilde a\) is \(C^1\) at \(x=1\), is positive, and is
constant for \(x\ge1+\delta\). If \(q\le0\), the required upper
differential inequality is immediate on \((1,2)\); if \(q>0\), then
\[
 \frac{x\widetilde a'(x)}{\widetilde a(x)}
 =xq\rho'(x-1)\le(1+\delta)q<\widetilde K_a
\]
where the derivative is nonzero. The additional monotonicity condition in
the SD case is local near \(x=0\) and is unchanged. Finally, rescaling
\(x=2\xi\) reduces~\eqref{sys:extended} to an interval of unit length and
preserves these properties up to fixed structural constants.

Since \(\widetilde c=c\) on \((0,1)\) and \(\widetilde c=0\) on
\((1,2)\), one has
\[
 \|\widetilde c(t,\cdot)\|_
 {L^\infty((0,2)\setminus\overline{\omega_e})}
 =\|c(t,\cdot)\|_{L^\infty(0,1)}.
\]
The asserted decay is therefore precisely~\eqref{cond:decay-bd}.
\end{proof}

\begin{proof}[Proof of Theorem~\ref{thm:boundary}]
Apply Theorem~\ref{thm:internal}, after the fixed spatial rescaling described
above, to system~\eqref{sys:extended}. There exists
\(\widetilde u\in L^2((0,T)\times\omega_e)\) such that the corresponding
solution satisfies
\[
 \widetilde y(T,\cdot)=0\qquad\text{in }(0,2).
\]
Define
\[
 y=\widetilde y|_{(0,T)\times(0,1)},\qquad
 v(t)=\widetilde y(t,1).
\]
The equation is uniformly parabolic on a fixed neighborhood
\(I_\varepsilon=(1-\varepsilon,1+\varepsilon)\) of \(x=1\).
Consequently, the one-dimensional trace inequality and the energy estimate
for~\eqref{sys:extended} give
\begin{align*}
 \|v\|_{L^2(0,T)}^2
 &\le C_\varepsilon\int_0^T
 \|\widetilde y(t,\cdot)\|_{H^1(I_\varepsilon)}^2\,dt\\
 &\le C\int_0^T\int_{I_\varepsilon}
 \bigl(|\widetilde y|^2+
       \widetilde a|\widetilde y_x|^2\bigr)\,dx\,dt<\infty.
\end{align*}
Thus \(v\in L^2(0,T)\). By restriction, \(y\) satisfies
\eqref{sys:boundary} in the variational sense, with initial datum \(y_0\),
history \(\Theta\), and boundary value \(y(t,1)=v(t)\). Finally,
\[
 y(T,\cdot)=\widetilde y(T,\cdot)|_{(0,1)}=0.
\]
The bounded dependence of \(\widetilde u\) on
\((\widetilde y_0,\widetilde\Theta)\), together with the energy and trace
estimates above, yields
\[
 \|v\|_{L^2(0,T)}
 \le C_T\|(y_0,\Theta)\|_{M_\tau}.
\]
This proves the theorem.
\end{proof}

\section{Conclusion}
We have proved null controllability for a one-dimensional degenerate
parabolic equation with a time-varying delay. The control may act either on
an interior subinterval or through a Dirichlet condition at the nondegenerate
endpoint. For the distributed problem, a degenerate Carleman estimate applies
on the final interval where the advanced term vanishes. A moving-limit energy,
adapted to the inverse delay map, propagates this estimate and leads to the
required range inclusion. Extending the equation beyond the controlled
endpoint reduces the boundary problem to the distributed one. The terminal
decay condition on the delayed coefficient is needed because reaching zero
before the terminal time does not generally keep a delayed state at rest.
Taking $\tau(t)\equiv h$ recovers the constant-delay equation.


\begin{thebibliography}{AKBDM14}
\bibitem[AKBDM14]{A5} F. Ammar-Khodja, C. Bouzidi, C. Dupaix, L. Maniar, \emph{Null controllability of retarded parabolic equations}, Math. Control Relat. Fields, 4 (2014), 1--15.
\bibitem[AF21]{AllalFragnelli} B. Allal, G. Fragnelli,
\emph{Null controllability of degenerate parabolic equation with memory},
Math. Methods Appl. Sci., 44 (2021), 9163--9190.
\bibitem[AFS22]{AllalFragnelliSalhi} B. Allal, G. Fragnelli, J. Salhi,
\emph{Controllability for degenerate/singular parabolic systems involving
memory terms}, Discrete Contin. Dyn. Syst. Ser. S, 15 (2022), 3445--3480.
\bibitem[Art69]{Artola} M. Artola, \emph{Sur les perturbations des \'equations d'\'evolution: Application \`a des probl\`emes avec retard}, Ann. Sci. \'Ecole Norm. Sup. (4), 2 (1969), 137--253.
\bibitem[ACF06]{BOU} F. Alabau-Boussouira, P. Cannarsa, G. Fragnelli, \emph{Carleman estimates for degenerate parabolic operators with applications to null controllability}, J. Evol. Equ., 6 (2006), 161--204.
\bibitem[CZ95]{Curt&Zwart} R. F. Curtain, H. Zwart, \emph{An Introduction to Infinite-Dimensional Linear Systems Theory}, Texts in Applied Mathematics 21, Springer, 1995.
\bibitem[AFM19]{benfama} E. M. Ait Ben Hassi, M. Fadili, L. Maniar, \emph{On algebraic conditions for null controllability of some coupled degenerate systems}, Math. Control Relat. Fields, 9 (2019), 77--95.
\bibitem[FM17]{FadiliManiar} M. Fadili, L. Maniar, \emph{Null controllability of $n$-coupled degenerate parabolic systems with $m$-controls}, J. Evol. Equ., 17 (2017), 1311--1340.
\bibitem[Nak86]{Nakagiri} S.-I. Nakagiri, \emph{Optimal control of linear retarded systems in Banach spaces}, J. Math. Anal. Appl., 120 (1986), 169--210.
\bibitem[NY89]{Nakagiri&Yamamoto} S.-I. Nakagiri, M. Yamamoto, \emph{Controllability and observability of linear retarded systems in Banach spaces}, Int. J. Control, 49 (1989), 1489--1504.
\bibitem[Maj23]{Majumdar} S. Majumdar,
\emph{Asymptotic behavior of the linearized compressible barotropic
Navier--Stokes system with a time varying delay term in the boundary or
internal feedback}, Math. Methods Appl. Sci., 46 (2023), 17288--17312.
\bibitem[WZ20]{WuZhang} H.-N. Wu, X.-W. Zhang,
\emph{Static output feedback stabilization for a linear parabolic PDE system
with time-varying delay via mobile collocated actuator/sensor pairs},
Automatica, 117 (2020), 108993.
\bibitem[ZZ18]{ZhouZhang} X. Zhou, M. Zhang,
\emph{On the controllability of a class of degenerate parabolic equations
with memory}, J. Dyn. Control Syst., 24 (2018), 577--591.
\bibitem[Zab95]{Zab} J. Zabczyk, \emph{Mathematical Control Theory}, Birkh\"auser, 1995.
\end{thebibliography}
\end{document}